\documentclass[12 pt]{article}

\usepackage{amsmath,amssymb,amsfonts,latexsym,float,graphics,epsfig}

\usepackage{setspace}
\usepackage[top=2cm, bottom=2cm, left=2cm, right=2cm]{geometry}
\usepackage[colorlinks]{hyperref}
\begin{document}

\renewcommand\theequation{\arabic{section}.\arabic{equation}}
\catcode`@=11 \@addtoreset{equation}{section}\newtheorem{axiom}{Definition}[section]
\newtheorem{theorem}{Theorem}[section]
\newtheorem{axiom2}{Example}[section]
\newtheorem{lem}{Lemma}[section]
\newtheorem{prop}{Proposition}[section]
\newtheorem{cor}{Corollary}[section]
\newtheorem{remark}{Remark}
\newcommand{\be}{\begin{equation}}
\newcommand{\ee}{\end{equation}}

\title{$3D$-flows Generated by the Curl of a Vector Potential \& Maurer-Cartan Equations} \maketitle
\begin{center}

O\u{g}ul Esen\footnote{E-mail: oesen@gtu.edu.tr}\\
Department of Mathematics, \\ Gebze Technical University, 41400 Gebze,
Kocaeli, Turkey.

\bigskip

Partha Guha\footnote{E-mail: partha.guha@ku.ac.ae}\\
Department of Mathematics,
Khalifa University\\
P.O. Box 127788, Zone -1 Abu Dhabi, UAE. 

\bigskip

Hasan G\"{u}mral\footnote{E-mail: hgumral@yeditepe.edu.tr}\\
Department of Mathematics, \\ 
Yeditepe University, 34755 Ata{\c{s}}ehir, {\.{I}}stanbul, Turkey.

\end{center}

\date{ }

\smallskip

\smallskip

\begin{center}

\end{center}

\bigskip

\begin{abstract}
\textit{We examine $3D$ flows $\mathbf{\dot{x}}=\mathbf{v}({\bf x})$ admitting vector identity $M\mathbf{v} = \nabla \times \mathbf{A}$ for a multiplier $M$ and a potential field $\mathbf{A}$. It is established that, for those systems, one can complete the vector field $\mathbf{v}$ into a basis fitting an $\mathfrak{sl}(2)$-algebra. Accordingly, in terms of covariant quantities, the structure equations determine a set of equations in Maurer-Cartan form. This realization permits one to obtain the potential field as well as to investigate the (bi-)Hamiltonian character of the system. The latter occurs if the system  has a time-independent first integral. In order to exhibit the theoretical results on some concrete cases, three examples are provided, namely the Gulliot system, a system with a non-integrable potential, and the Darboux-Halphen system in symmetric polynomials. 
}
\end{abstract}

\tableofcontents
\setlength{\parskip}{4mm}

\onehalfspacing

\section{Introduction}

In a recent work \cite{EsGu20}, we have presented (bi)Hamiltonian analysis of $3D$ dynamical systems $\mathbf{\dot{x}}=\mathbf{v}({\bf x})$ where the velocity field is the curl of 
a vector  potential that is $\mathbf{v} = \nabla \times \mathbf{A}$ for some $\mathbf{A}$. The analysis was divided into two cases according to the (Frobenius) integrability of the potential vector field $\mathbf{A}$. Two examples have been provided; a bi-Hamiltonian system admitting a non-integrable potential and a non-Hamiltonian system admitting an integrable potential. If a system possesses bi-Hamiltonian character then the flow is the line of intersection of two surfaces determined by  the Hamiltonian functions \cite{Barbarosie,CGR,Morando}. This geometric realization is a particular instance of the  superintegrability \cite{Ev}.

In the present paper, we are addressing the same problem from an algebraic point of view. The goal is to determine a potential field $\mathbf{A}$ for a given system $\mathbf{\dot{x}}=\mathbf{v}({\bf x})$ satisfying $\mathbf{v} = \nabla \times \mathbf{A}$. Evidently, this holds for volume preserving flows. For the other case, the curl identity needs to be upgraded to $M\mathbf{v} = \nabla \times \mathbf{A}$ where $M$ is a conformal factor called Jacobi last multiplier, see \cite{ChGuKh09,Nu}.
More concretely, in this work, we shall argue that the existence of a vector potential is manifesting a representation of $\mathfrak{sl}(2)$ algebra of vector fields spanning the configuration space. Referring to the  dualization between a vector field and a one-form section, we carry this algebra to the level of differential forms. This leads us to determine Maurer-Cartan type equations closing the exterior algebra of sections. We call this geometrization as Maurer-Cartan $\mathfrak{sl}(2)$ algebra of curl vector fields. Interestingly, in the Maurer-Cartan $\mathfrak{sl}(2)$ algebra, one of the one-form sections is not necessarily integrable. This permits us to claim that the present analysis is applicable even for non-integrable cases. Further, referring to the conformal invariance of the algebra, we shall state that a perturbation is possible taking a non-integrable potential to an integrable one. 

The paper is organized into $5$ main sections. In the following one, we shall state some necessary background of $3D$ systems. Then, in Section \ref{Sec-Curl}, we shall focus on  $3D$ systems admitting vector potentials. In Section \ref{Sec-SE}, we shall  provide    Maurer-Cartan $\mathfrak{sl}(2)$ algebra of one-forms. Referring to the classical duality, in Section \ref{Sec-DS}, it will be shown that $3D$ systems admitting vector potentials are inducing the $\mathfrak{sl}(2)$ algebra of vector fields. Three examples will be provided in Section \ref{Sec-IL} including the Gulliot system, a system with a non-integrable potential, and the Darboux-Halphen system in symmetric polynomials. 

\section{Bi-Hamiltonian Dynamics in $3D$}

Let $(\mathcal{P},\{\bullet,\bullet\})$ be a $3$-dimensional Poisson manifold equipped with a Poisson bracket $\{\bullet,\bullet\}$. The Hamilton's equation generated by a Hamiltonian function $H$ is defined to be 
\begin{equation} \label{HamEq1}
\mathbf{\dot{x}}=\{\mathbf{x},H\},
\end{equation}
for local coordinates $(\mathbf{x})$ on $\mathcal{P}$. 
In $3$-dimensions, we can replace the role of a Poisson bracket with a Poisson vector $\mathbf{J}$, \cite{EsGhGu16, Gum1, GuNu93}. In this case, the Jacobi identity turns out to be the following vector equation
\begin{equation} \label{jcbv}
\mathbf{J}\cdot(\nabla\times\mathbf{J})=0,
\end{equation}
whereas the Hamilton's equation (\ref{HamEq1}) takes
the particular form 
\begin{equation} \label{HamEq3}
\mathbf{\dot{x}}=\mathbf{J}\times\nabla H.
\end{equation}
Here, $\nabla H$ is the gradient of $H$.
The following theorem is exhibiting all possible solutions of the Jacobi identity (\ref{jcbv}) so it characterizes Poisson structures in $3$-dimensions, \cite{HB1,HB2,HB3}.
\begin{theorem} \label{ss}
The general solution of the vector equation (\ref{jcbv}) is 
$\mathbf{J}= \left({1}/M \right) \nabla F$ for arbitrary functions $M$
and $F$. \end{theorem}

The existence of scalar multiple ${1}/M$ in the solution is a manifestation of  conformal invariance of the identity (\ref{jcbv}). In the literature, $M$ is called Jacobi's last multiplier \cite{Jac1,Jac2}. In this picture, a Hamiltonian system has the following generic form
\begin{equation} \label{HamEq4}
\mathbf{\dot{x}}= \frac{1}{M} \nabla F \times \nabla H.
\end{equation}

A dynamical system is bi-Hamiltonian if it admits two different Hamiltonian
structures
\begin{equation}
\mathbf{\dot{x}}= \{\mathbf{x},H_2\}_{1}=\{\mathbf{x},H_1\}_{2},\label{biHam}%
\end{equation}
with the requirement that the Poisson brackets $\{\bullet,\bullet \}_{1}$ and $\{\bullet,\bullet \}_{2}$ be
compatible \cite{Fe94,OLV}. That is any linear pencil $\{\bullet ,\bullet\}_{1}+c\{\bullet ,\bullet\}_{2}$ must satisfy the Jacobi
identity \cite{MaMo84,OLV}. In three dimensions,  a bi-Hamiltonian system can be put into the form
\begin{equation}
M\mathbf{\dot{x}}=\mathbf{J}_{1}\times \nabla H_{2}=\mathbf{J}_{2}\times
\nabla H_{1}.  \label{bi-Ham-}
\end{equation}%

Referring to the system (\ref{HamEq4}), we conclude that a Hamiltonian system in the form of (\ref{HamEq4}) is bi-Hamiltonian
\begin{equation} \label{bi-Ham}
M \mathbf{\dot{x}}=  \nabla H_1\times\nabla H_2=\mathbf{J}_{1}\times\nabla H_1=\mathbf{J}_{2}\times\nabla
H_2,
\end{equation}
where, the first Poisson vector $\mathbf{J}_{1}=-\nabla H_2$ whereas the second Poisson vector $\mathbf{J}_{2}= \nabla H_1$. The following theorem determines the Hamiltonian picture of three dimensional dynamical systems admitting an integral invariant. For the proof, we refer \cite{EsGhGu16,Gao}.
\begin{theorem} \label{ss2}
A three dimensional dynamical system $\mathbf{\dot{x}}=\mathbf{v}(\bf{x})$ having a
time independent first integral is bi-Hamiltonian if and only if there exist
a conformal factor, called Jacobi's last multiplier, $M$ which makes $M\mathbf{v}$ divergence-free.
\end{theorem}

\section{Curl Fields}\label{Sec-Curl}

We have depicted the generic form of the $3D$ dynamical systems in \eqref{HamEq4}. In this section, for a given $3D$ dynamical system $\mathbf{v}=\mathbf{\dot{x}}$, we examine the existence of a potential vector field $\mathbf{A}$ determining the dynamics up to some conformal factor.

Assume that, the dynamical field is in the bi-Hamiltonian  form $\mathbf{v}=\nabla
H_{1}\times \nabla H_{2}$ where the multiplier $M$ is being unity. In this case, we can write $\mathbf{v}$ as a
curl vector $\mathbf{v}=\nabla \times \mathbf{A}$ with $\mathbf{A}%
=H_{1}\nabla H_{2}$ so that $\nabla \cdot \mathbf{v}=0$. To see this realization in a covariant formulation, consider the standard coordinates $(x,y,z)$ on the space, and define a two-form
\begin{equation}
\mathbf{v}\cdot d%
\mathbf{x}\wedge d\mathbf{x}:=
\iota _{v}\left( dx\wedge dy\wedge dz\right),
\end{equation}
where $\iota_{v}$ is the interior product by the vector field $v=\mathbf{v}\cdot \nabla$. The volume preserving character of the field $\mathbf{v}$ can be recorded in terms of the Lie derivative as
\begin{equation}
\mathcal{L}_{v}(dx\wedge dy\wedge dz)= d\iota _{v}(dx\wedge dy\wedge dz)=\left( \nabla \cdot \mathbf{v}%
\right) dx\wedge dy\wedge dz=0,
\end{equation}
where we have employed the Cartan's identity $\mathcal{L}_{v}=d\iota _{v}+ \iota _{v}d$. If $\mathbf{v}$ is not divergence free, that is  
 $\nabla \cdot \mathbf{v}\neq 0$, then assume an
invariant volume 
\begin{equation}
\ast 1=M dx\wedge dy\wedge dz
\end{equation}
involving a conformal factor $M$. Here, $\ast$ is the Hodge star operator.  This is a manifestation of Theorem \ref{ss2}. In this case, one can recast the conservation of the invariant volume as 
\begin{equation}
\mathcal{L}_{v}\left( \ast 1\right)= d\iota _{v}\left( \ast 1\right) =\left( \nabla \cdot M \mathbf{v}%
\right) dx\wedge dy\wedge dz. 
\end{equation}
Accordingly, a vector potential in form $\mathbf{A}=H_{1}\nabla H_{2}$ satisfies the equation $
M  \mathbf{v}=\nabla \times \mathbf{A}$.

To switch to a covariant picture, we define the Poisson one-form
\begin{equation}
\mathbf{J} \longrightarrow J=J_{i}dx^{i}.
\end{equation}
In this case, the Jacobi identity is given by the Frobenius integrability  condition $J\wedge dJ=0$. We have coordinate independent manifestation of the dynamics
\begin{equation}
\iota _{v}\left( \ast 1\right) =J^{(1)}\wedge J^{(2)}
\end{equation}
where the left hand side can be written as
\begin{eqnarray*}
\iota _{v}\left( \ast 1\right)  &=&\left(M \mathbf{v}\right)
\cdot d\mathbf{x}\wedge d\mathbf{x=}\nabla \times \mathbf{A}\cdot d\mathbf{x}%
\wedge d\mathbf{x} \\
&=&d\left( \mathbf{A}\cdot d\mathbf{x}\right) =d \gamma
\end{eqnarray*}%
for a one-form (potential) $\gamma=\mathbf{A}\cdot d\mathbf{x}$ of $M\mathbf{v}$. Thus, casting $M\mathbf{v}$ into bi-Hamiltonian form is the same as writing a Maurer-Cartan like equation%
\begin{equation}
d\gamma=J^{(1)}\wedge J^{(2)}.
\end{equation}
 
\section{Structure Equations}\label{Sec-SE}

Given a three dimensional dynamical system $\mathbf{\dot{x}}=\mathbf{v}$, our interest is investigating a potential one-form for $\mathbf{v}$. We first assume the existence of an integrable one-form potential. Then, referring to this one-form, we shall construct an algebra in the space of one-form sections on $\mathbb{R}^3$. Later, the case of non-integrable potential one-forms will be examined in the light of the present discussion.    

Assume that $\gamma$ represents an integrable  potential one-form for $\mathbf{v}$. The integrability condition $\gamma \wedge d\gamma =0$ implies that there exists a
one-form $\alpha $, mimicking the role of an integrating factor, such that 
\begin{equation}\label{111}
d\gamma =2\alpha \wedge \gamma .
\end{equation}
Taking the exterior derivative of \eqref{111}, we arrive at the following
\begin{equation}
2d\alpha \wedge \gamma =2\alpha \wedge 2\alpha \wedge \gamma =0.
\end{equation}
This identity determines two possibilities. First, $\alpha$ is a closed one-form, then we can integrate and we are done. Second, $d\alpha\neq 0$ and we have
\begin{equation}\label{222}
d\alpha =\gamma \wedge \beta
\end{equation} 
for some one-form $\beta $. Since, we have assumed that $\alpha$ is not closed, $\gamma$ and $\beta$ are linearly independent. We further assume that the set $\{\alpha,\beta,\gamma\}$ determines a basis for the one-form sections.  An implication of this is $\alpha \wedge d\alpha\neq 0$ which says that $\alpha $ is not integrable. We shall comment on this after conformal invariance of the structure equations is obtained. From the exterior derivative of \eqref{222}
\begin{equation}\label{222-}
2\alpha \wedge \gamma \wedge \beta -\gamma \wedge d\beta =0.
\end{equation}
and linear independence of basis one-forms, we obtain 
\begin{equation}\label{333}
d\beta =-2\alpha \wedge \beta .
\end{equation}
Thus, by starting with an integrable one-form $\gamma$, we obtained a linearly independent basis satisfying the structure equations \eqref{111}, \eqref{222} and \eqref{333}. We collect these Maurer-Cartan type equations  in the following theorem while exhibiting $\mathfrak{sl}(2)$ algebra character of the system.  
\begin{theorem}\label{Basis}
An integrable one-form $\gamma$ determines $3$-dimensional basis satisfying 
\begin{equation}\label{MC-alg}
d\beta =-2\alpha \wedge \beta ,\qquad d\alpha =\gamma \wedge \beta ,%
\qquad d\gamma =2\alpha \wedge \gamma,
\end{equation}
where $d\alpha \neq 0$,
\end{theorem}
The algebra given in \eqref{MC-alg} admits a symmetry by being invariant under some conformal transformations. 
\begin{theorem}\label{MC-inv}
The Maurer-Cartan  system \eqref{MC-alg} is invariant under conformal transformations
\begin{equation}
\beta \mapsto \rho \beta, \qquad \alpha \mapsto \alpha -\frac{1 }{2 }d\ln \rho, \qquad \gamma \mapsto \frac{\gamma }{\rho },
\end{equation}
for a non-vanishing function $\rho $. 
\end{theorem}
To prove this assertion, we start with $\beta \mapsto \rho \beta $ and compute  
\begin{equation*}
d\left( \rho \beta \right)  = d\rho \wedge \beta +\rho \wedge d\beta  = \frac{d\rho }{\rho }\wedge \left( \rho \beta \right) +\rho \left(
-2\alpha \wedge \beta \right)  = -\big( 2\alpha -\frac{d\rho }{\rho }\big) \wedge \left( \rho \beta
\right) 
\end{equation*}
which defines $\alpha \mapsto \alpha -\frac{d\rho }{2\rho }$. Differentiating  
\begin{equation*}
d\big( \alpha -\frac{d\rho }{2\rho }\big)  = d\alpha =\gamma \wedge
\beta =\gamma \wedge \frac{\rho }{\rho }\beta = \frac{\gamma }{\rho }\wedge \left( \rho \beta \right)
\end{equation*}
gives the form of transformation $\gamma \mapsto \frac{\gamma }{\rho }$. Further differentiation closes the algebra
\begin{equation*}
d\big( \frac{\gamma }{\rho }\big) =\frac{d\gamma }{\rho }-\frac{1}{%
\rho ^{2}}d\rho \wedge \gamma =\frac{2\alpha \wedge \gamma }{\rho }-\frac{%
d\rho }{\rho }\wedge \frac{\gamma }{\rho } = 2\big( \alpha -\frac{1}{2}\frac{d\rho }{\rho }\big) \wedge \frac{%
\gamma }{\rho }.
\end{equation*}

\textbf{Non-integrable Integrating Factor.} In the structure equations, the one-form $\alpha$ appears as integrating factor for integrable one-forms $\beta$ and $\gamma$. Yet $\alpha$ itself is non-integrable, because
\begin{equation*}
\alpha\wedge d\alpha=\alpha\wedge\gamma\wedge\beta\neq 0
\end{equation*}
for an orientable three manifold. As an integrating factor, we can seek for functions $f$ and $g$ which will make $\alpha$ integrable in 
\begin{equation}
\begin{split}
d\beta & = -2(\alpha+f\beta)\wedge \beta = -2\alpha\wedge\beta
\\
d\gamma & = 2(\alpha + g \gamma)\wedge \gamma = 2\alpha\wedge\gamma.
\end{split}
\end{equation}
that is we require $f$, $g$ to satisfy  
\begin{equation}
\begin{split}
(\alpha+f\beta) \wedge d(\alpha+f\beta)&=0
\\
(\alpha + g \gamma) \wedge d(\alpha + g \gamma) & =0. 
\end{split}
\end{equation}
These conditions imply linear first order PDEs for the functions $f$ and $g$
\begin{equation}
\begin{split}
\alpha  \wedge  (\gamma + df) \wedge \beta   &=0
\\
\alpha  \wedge  (-\beta + dg) \wedge \gamma   &=0
\end{split}
\end{equation}
which can always locally solvable. This means that, in the decomposition of two-forms $d\beta$ and $d\gamma$ we can replace non-integrable integrating factor $\alpha$ with integrable ones $\alpha+f\beta$ in $d\beta$ and $\alpha+g\gamma$ in $d\gamma$ to make them integrable and hence Poisson one-forms for vector fields corresponding to $d\beta$ and $d\gamma$. 

\textbf{Non-integrable Potential One-form.} Here, we assume that $\mathbf{v}$ admits a non-integrable potential vector. In this case, we simply can identify the potential one-form to be $\alpha$ which already presents in $\mathfrak{sl}(2)$-structure. Then, the locally bi-Hamiltonian form of $\mathbf{v}$ is one of the Maurer-Cartan equations
\begin{equation*}
\iota_v(dx \wedge dy\wedge dz)=d\alpha =\gamma\wedge \beta
\end{equation*}
with the potential $\alpha$ being non-integrable 
\begin{equation*}
\alpha \wedge d\alpha =\alpha \wedge \gamma\wedge \beta \neq 0.
\end{equation*}

\section{Dynamical System}\label{Sec-DS}

Suppose the system  $\mathbf{\dot{x}=v}$ comes along with  $
\mathbf{u}$ and $\mathbf{w}$ constituting an $\mathfrak{sl}(2)$ algebra
\begin{equation}\label{sl-2-v}
[u,v]=2v,\qquad [u,w]=-2w, \qquad [v,w]=u
\end{equation}
where $u=\mathbf{u\cdot \nabla }$, $v=\mathbf{v\cdot \nabla }$ and $w=%
\mathbf{w\cdot \nabla }$. Here, the brackets are the Jacobi-Lie bracket of vector fields. In this realization, the invariant volume density is
\begin{equation} \label{LM}
\frac{1}{M}=\left( \mathbf{v\times u}\right) \cdot \mathbf{w}.
\end{equation}
We define the dual one-forms
\begin{equation}\label{a-b-c}
\alpha =M  ( \mathbf{w\times v} ) \cdot d\mathbf{x},\qquad
\beta = M ( \mathbf{u\times w} ) \cdot d\mathbf{x},\qquad
\gamma =M  ( \mathbf{v\times u} ) \cdot d\mathbf{x}
\end{equation}
which can easily be shown to  satisfy
\begin{equation}
\iota _{v}\beta =\iota _{u}\alpha =\iota _{w
}\gamma =1,
\end{equation}
and all the other possible couplings are identically zero. In order to see the correspondence between $\mathfrak{sl}(2)$ algebra in \eqref{sl-2-v} and the one exhibited in \eqref{MC-alg}, it is enough to consider the very definition of the exterior derivative. For a one-form $\omega$, it is given by
\begin{equation}
d\omega  ( u,v ) =u( \iota _{v}\omega ) -v ( \iota
_{u}\omega  ) -\iota _{\lbrack u,v]}\omega. 
\end{equation}

To show that $\mathfrak{sl}(2)$ is the natural structure regarding $\mathbf{v}$ as a curl vector field, we consider the one-form $\gamma $ in \eqref{a-b-c}. To show that, $d\gamma $ is a curl
expression for $\mathbf{v}$, we compute
\begin{equation} \label{All-1}
\begin{split}
\nabla \times \left( M \mathbf{v\times u}\right)  &= \nabla M \times
\left( \mathbf{v\times u}\right) +M \nabla \times \left( \mathbf{v\times u
}\right)  \\
&= (\mathbf{u} \cdot\nabla M ) \mathbf{v}-(\mathbf{v} \cdot\nabla M )    \mathbf{u}
+M \left( \nabla \cdot \mathbf{u}\right) \mathbf{v}-M  \left( \nabla \cdot 
\mathbf{v}\right) \mathbf{u}+M \overrightarrow{[u,v]} \\
&= \nabla \cdot \left( M  \mathbf{u}\right) \mathbf{v}-\nabla \cdot \left(
M \mathbf{v}\right) \mathbf{u}+2 M \mathbf{v} = 2M \mathbf{v}
\end{split}
\end{equation}
 due to the invariance of $M$. Similarly, one can compute that%
\begin{equation} \label{All-2}
\nabla \times \left( M  \mathbf{u\times w}\right) =2M\mathbf{w}, \qquad 
\nabla \times \left( M  \mathbf{v\times w}\right) =-M\mathbf{u}.
\end{equation}
Recall that for a one-form $\mathbf{A}\cdot d\mathbf{x}$ we have $d ( \mathbf{A}\cdot d\mathbf{x} ) =\left( \nabla \times \mathbf{A}\right) \cdot d\mathbf{x}\wedge d\mathbf{x}$. So that, the left-hand sides of the Maurer-Cartan system in \eqref{MC-alg} are the curls whereas  the right-hand sides
decompose these curls into two potentials that are Poisson vectors. \eqref{All-1} and \eqref{All-2} are coefficients of $d\mathbf{x} \wedge d\mathbf{x}$ in $d\gamma$, $d\beta$ and $d\alpha$, respectively. To verify the right hand sides of Maurer-cartan equations, we take as an example $d\gamma =2\alpha \wedge \gamma $ and compute 
\begin{equation} 
\begin{split}
2M d\mathbf{x} \wedge d\mathbf{x}&=2 M( \mathbf{w\times v} ) \cdot d\mathbf{x}  \wedge M ( \mathbf{v\times u} )\cdot d\mathbf{x}
\\
&= 2M^2  ( \mathbf{w\times v} ) \times ( \mathbf{v\times u} )\cdot d\mathbf{x} \wedge d\mathbf{x}
\\
&= 2M^2 \mathbf{v} (\mathbf{u}\cdot \mathbf{w\times v}) \cdot 
d\mathbf{x} \wedge d\mathbf{x}
=2M \mathbf{v} \cdot d\mathbf{x} \wedge d\mathbf{x}
\end{split}
\end{equation}
by definition of the multiplier $M$.

\textbf{Heisenberg Algebra.} We consider a basis $\{\omega^1,\omega^2,\omega^3\}$ for the space of one-form sections and  assume that the following structure equations hold
\begin{equation}\label{Heisenberg}
d\omega ^{1}=d\omega ^{3}=0,\qquad d\omega ^{2}=\omega ^{3}\wedge
\omega ^{1}.
\end{equation}
Referring to the standard coordinates $(x,y,z)$, we can write the local realizations 
\begin{equation}\label{omega}
\omega ^{1}=dx, \qquad \omega ^{2}=dy-xdz, \qquad \omega ^{3}=dz.
\end{equation}
For each one form $\omega=\boldsymbol{\omega}\cdot d\mathbf{x} $, we associate a covector $\boldsymbol{\omega}$. Then referring to \eqref{omega}, we arrive at the following set
\begin{equation}
\boldsymbol{\omega}^{1}= ( 1,0,0 ) , \qquad \boldsymbol{\omega}
^{2}=\left( 0,1,-x\right) ,\qquad \boldsymbol{\omega}^{3}=\left(
0,0,1\right)
\end{equation}
satisfying $M =\boldsymbol{\omega}^{1}\times \boldsymbol{\omega}^{2}\cdot \boldsymbol{\omega}^{3}=1$. Accordingly, we define the following basis 
\begin{equation}
u  =\boldsymbol{\omega}^{1}\times \boldsymbol{\omega}^{2}\cdot \nabla =x%
\frac{\partial }{\partial y}+\frac{\partial }{\partial z} 
, \qquad 
v  =\boldsymbol{\omega}^{3}\times \boldsymbol{\omega}^{1}\cdot \nabla =%
\frac{\partial }{\partial y} , \qquad 
w  = \boldsymbol{\omega}^{2}\times \boldsymbol{\omega}^{3}\cdot \nabla = \frac{\partial }{\partial x} 
\end{equation}%
satisfying
\begin{equation}\label{alg-Hei}
\lbrack v,w]=[v,u]=0, \qquad [w,u]=v.
\end{equation}%
If $\mathbf{v}$ is the given dynamical system, these equations characterize $\mathbf{v}$ as having two symmetries whose commutator produces the dynamics.
Note the dualities%
\[
\iota _{w}\mathbf{\omega }^{1}=\iota _{v}\mathbf{\omega }%
^{2}=\iota _{u}\mathbf{\omega }^{3}=1.
\]%
It follows from \eqref{Heisenberg}  that $\boldsymbol{\omega}^{1}$ and $\boldsymbol{\omega}^{3}$ are
conserved covariants and $\mathbf{v}$, expressed as a curl, is bi-Hamiltonian with these invariants. 
This
shows that well-known bi-Hamiltonian systems with two integrals of motion can
in fact be manifested in Heisenberg algebra \eqref{Heisenberg}. 
These local coordinates are indeed the case of final quadrature after having conserved quantities, 
say $H_1$ and $H_2$, and reducing the system by eliminating two coordinates. That is
\[
\mathbf{x} =\left( H_1,y,H_2\right), \qquad  \mathbf{\dot{x}}%
=\left( 0,\dot{y}(y,H_1,H_2),0\right).
\]
\section{Examples} \label{Sec-IL}

\subsection{Guillot System}
As a first example, we consider the Guillot System \cite{Gu}
\begin{equation} \label{Gu}
\dot{x} =x^{2}+y^{4},  \qquad
\dot{y} =xy, \qquad 
\dot{z} = 2y^{2}z-xz .
\end{equation}
In order to investigate a potential vector field for this system, the first step is to 
determine the vector field $v$ generating the system \eqref{Gu} and then complete it to a $3$ dimensional basis satisfying  structure equations for $\mathfrak{sl}(2)$ that is \eqref{sl-2-v} which has already been done by Guillot in search of vector fields of Darboux-Halphen type 
\begin{equation}
\begin{split}
v &=  ( x^{2}+y^{4} ) \frac{\partial }{\partial x}+xy\frac{\partial 
}{\partial y}+  ( 2y^{2}z-xz ) \frac{\partial }{\partial z}, \\
u &= 2x\frac{\partial }{\partial x}+y\frac{\partial }{\partial y}-z\frac{
\partial }{\partial z} ,\\
w &= -\frac{\partial }{\partial x}.
 \end{split}
\end{equation}
The reciprocal of the multiplier is given in 
\eqref{LM}. For the present system, it is computed to be 
\begin{equation}\label{M-Gu}
M = \frac{1}{ \mathbf{v\times u}  \cdot \mathbf{w} }=\frac{1}{2zy^{3}}
\end{equation}
provided that $2zy^{3}$ never vanishes. It is now straight forward to check that $Mv$ is a divergence free vector field. Referring to the identifications given in \eqref{a-b-c}, we compute the following one-form sections
\begin{equation}\label{a-b-c-G}
\begin{split}
\alpha &= \frac{2y^2-x}{2y^3} dy-\frac{x}{2zy^{2}}dz,\\
\beta &= \frac{1}{2y^3} dy + \frac{1}{2zy^{2}}dz, \\
\gamma &= -dx + \frac{x^2-y^4-4xy^2}{2y^3} dy + \frac{y^4-x^2y}{2zy^{2}}dz,
 \end{split}
\end{equation}
respectively. It is now the matter of a direct calculation to verify that the forms in \eqref{a-b-c-G} are satisfying the $\mathfrak{sl}(2)$  Maurer-Cartan equations \eqref{MC-alg} stated in Theorem \ref{Basis}. See also that if we write $\gamma=\mathbf{A}\cdot d\mathbf{x}$ then it is possible now to establish that $d\gamma=M\mathbf{v}\cdot d\mathbf{x}\wedge d\mathbf{x}$ which reads the vector $M\mathbf{v}$ as the curl of the potential field
\begin{equation}
\mathbf{A}=(-1,\frac{x^2-y^4-4xy^2}{2y^3} ,\frac{y^4-x^2y}{2zy^{2}}). 
\end{equation}

It is immediate to see that the following function
\begin{equation} \label{H_1-Gu}
H_1=\frac{x^{2}}{y^{2}}-y^{2}
\end{equation}
is a first integral of the Guillot System \eqref{Gu} for $y\neq 0$. Along with the existence of the multiplier $M$ in \eqref{M-Gu}, in the light of Theorem \ref{ss2}, we argue that the Guillot System \eqref{Gu} admits two conserved quantities so that one can recast it in the bi-Hamiltonian form. This implies that, one may examine the Guillot System \eqref{Gu} in the framework of Heisenberg algebra \eqref{alg-Hei} as well. To have this, we compute another first integral
\begin{equation}\label{H_2-Gu}
H_2=\epsilon  \log \big ( \frac{\epsilon x+y^{2}}{y (
yz) ^{\epsilon /2}}  \big) 
\end{equation}
where $\log$ is the natural logarithm and  $\epsilon $ stands for $\pm 1$. A direct computation exhibits bi-Hamiltonian character of the dynamics 
 \begin{equation}\label{ext2}
 2\iota_{v} (dx\wedge dy \wedge dz) = \frac{1}{M} dH_2 \wedge dH_1
 \end{equation}
where $M$ is the Jacobi's last multiplier in \eqref{M-Gu}. Here the factor $2$ on the right hand side of \eqref{ext2} is a manifestation of the multiple $2$ on the right hand side of the first bracket in the $\mathfrak{sl}(2)$ algebra \eqref{sl-2-v}. The flow is the line of intersection of two surfaces determined by  the Hamiltonian functions $H_1$ in \eqref{H_1-Gu} and $H_2$ in \eqref{H_2-Gu}.

\subsection{Non-integrable Potential} 

As a simple example, we take the vector field $\mathbf{u}$ in Guillot system. Dual one-form is $\alpha$ which is non-integrable. This means that 
\begin{equation} 
-M \mathbf{u}  \cdot d\mathbf{x} \wedge d\mathbf{x} =d\alpha =\gamma \wedge \beta
 \end{equation}
admits a potential one-form
\begin{equation} 
\alpha=\frac{2y^2-x}{2y^3}dy - \frac{x}{2zy^2}dz = \boldsymbol{\alpha}\cdot d\mathbf{x} 
 \end{equation}
which is not integrable $\alpha\wedge d \alpha \neq 0$. Indeed, we have $M \mathbf{u}=\nabla \times \boldsymbol{\alpha}$ but
\begin{equation} 
\boldsymbol{\alpha} \cdot \nabla \times \boldsymbol{\alpha}
=\frac{1}{2zy^3} = M
 \end{equation}
as expected. Non-integrable potential leads to the emergence of cohomological element called Godbillon-Vey class. This is the three-form 
\begin{equation} 
\alpha\wedge d \alpha = ( \boldsymbol{\alpha} \cdot \nabla \times \boldsymbol{\alpha} ) \, dx \wedge dy \wedge dz = *1
 \end{equation}
 which is obviously closed (we are in dimension three) but not exact, i.e. there is no two-form whose derivative is $ *1$. In other words, there is no vector field $\mathbf
{B}$ that solves the equation $\nabla \cdot \mathbf
{B} = \boldsymbol{\alpha} \cdot \nabla \times \boldsymbol{\alpha} $. 
 
\subsection{Darboux-Halphen System in Symmetric Polynomials }

We start with Darboux-Halphen system \cite{Da} given by 
\begin{equation}\label{DH-class}
\begin{split}
\dot{t}_1 = t_2 t_3-t_1 t_2-t_1 t_3,\\
\dot{t}_2 = t_1 t_3-t_3 t_2-t_1 t_2,\\
\dot{t}_3 = t_1 t_2-t_3 t_1-t_3 t_2.
\end{split}
\end{equation} 
See, for example, \cite{MoMoNiRoTo18} for more recent discussions on this system. 
We introduce a new set of dependent variables 
\begin{equation}
 x:=-2\left( t_1+t_2+t_3\right), \qquad y:=4\left( t_1 t_2+t_2 t_3+t_1 t_3\right), \qquad z:=-8t_1 t_2 t_3
\end{equation} 
which take Darboux-Halphen  system \eqref{DH-class} into the following form 
\begin{equation}\label{DH}
\dot{x}=\frac{1}{2}y, \qquad \dot{y}=3z, \qquad \dot{z}=2xz-\frac{1}{2}y^{2}.
\end{equation} 
We denote the dynamics governed by the equations \eqref{DH} by a vector field $v$. We complete the vector field $v$ into a basis $(v,u,w)$ satisfying  the $\mathfrak{sl}(2)$ algebra exhibited in \eqref{sl-2-v}. A direct computation gives this basis as
\begin{equation}
\begin{split}
v &= \frac{1}{2}y\frac{\partial }{\partial x}+3z\frac{\partial }{\partial y}+\big( 2xz-\frac{1}{2}y^{2}\big) \frac{\partial }{\partial z}, \\
u &= 2x\frac{\partial }{\partial x}+4y\frac{\partial }{\partial y}+6z\frac{\partial }{\partial z}, \\
w &= -6\frac{\partial }{\partial x}-4x\frac{\partial }{\partial y}-2y\frac{\partial }{\partial z}.
\end{split}
\end{equation} 
Referring to the identity \eqref{LM}, the volume is computed to be 
\begin{equation}\label{mult-DH}
\frac{1}{M}=72xyz-16y^{3}+4x^{2}y^{2}-16x^{3}z-108z^{2}.
\end{equation} 
According to \eqref{a-b-c}, we get the one form sections as
\begin{equation}\label{a-b-c-DH}
\begin{split}
\alpha  &=M \big( ( 2xy^{2}+6yz-8x^{2}z ) dx+ (
12xz-4y^{2} ) dy+ ( 2xy-18z ) dz\big), 
\\
\beta  &=4M \big(  ( 6xz-2y^{2} ) dx+ ( xy-9z )
dy+ ( 6y-2x^{2} ) dz\big),
\\
\gamma  &=M  \big( ( 18z^{2}-8xy+2y^{3} ) dx+ (
4x^{2}z-xy^{2}-3yz ) dy+ ( 2y^{2}-6xz ) dz\big),
\end{split}
\end{equation} 
where $M$ is the multiplier in \eqref{mult-DH}. 
It is now straight forward to check that these forms are satisfying the equations \eqref{MC-alg} presented in Theorem \ref{Basis}. Polynomial character of the multiplier $M$ permits us to assign dimensions $[x]=1$, $[y]=2$, and $[z]=3$. Thus $[M]=-6$. Accordingly, one computes the dimensions of the one-form sections in \eqref{a-b-c-DH} as $[\alpha ]=0$, $[\beta ]=-1$ and $[\gamma]=1$, respectively. This is obeying $\mathfrak{sl}(2)$ algebra  realization of the one-form sections. 
The exterior derivative $d\gamma $ will define the curl
potential for the dynamics $v$. Since Darboux-Halphen is known not to admit any polynomial 
integrals, this structure shows only that bi-Hamiltonian structure exists.

\end{document}